\documentclass[12pt,leqno]{article}
\usepackage{amsmath,amsthm}

\def\init{\setcounter{equation}{0}}
\newtheorem{theorem}{Theorem}[section]

\newcommand{\R}{{\bf R}}
\newcommand{\C}{{\bf C}}

\newtheorem{lemma}{Lemma}[section]

\newcommand{\e}{{\varepsilon}}

\title{On non-abelian Radon transform.
\author{G.Eskin, \ \ \  Department of Mathematics, UCLA,\\ Los Angeles,
CA 90095-1555, USA. \ E-mail: eskin@math.ucla.edu}
\footnote{This work was done when the author was participating in the Inverse 
Problems
program at IPAM, UCLA, Fall 2003.}
}

\begin{document}
\maketitle
\begin{abstract}
We consider the inverse prolem of the recovery of a gauge field in $\R^2$
modulo gauge transformations from the non-abelian Radon transform.  
A global uniqueness theorem is proven for the case when the gauge field
has a compact support.  Extensions to the attenuated non-abelian Radon 
transform in $\R^2$ and applications to the inverse scattering problem
for the Schr\"{o}dinger equation in $\R^2$ with non-compact Yang-Mills 
potentials are studied. 
\end{abstract}

\section{Introduction.}
\label{section 1}
\init

Let $A_j(x),\ 0\leq j\leq 2,$ be $C^\infty$  $m\times m$ matrice in
$\R^2$ with compact support:
$\mbox{supp\ } A_j(x)\subset B_R=\{x:|x|<R\},\ 0\leq j\leq 2.$
Denote $\theta=(\theta_1,\theta_2)\in S^1$.  Let $c_0(x,\theta)$ be
the matrix solution of the differential equation
\begin{equation}                           \label{eq:1.1}
\theta\cdot \frac{\partial c_0(x,\theta)}{\partial x}=
\left(A_1(x)\theta_1+A_2(x)\theta_2 +A_0(x)\right)c_0(x,\theta),
\end{equation}
such that
\[
c_0(x+s\theta,\theta)\rightarrow I_m\ \ \ \mbox{when\ }\ \ 
\ s\rightarrow -\infty.
\]

Denote by $S(A)$
the limit of $c(x+s\theta,\theta)$ when $s\rightarrow \infty$.  Note
that $S(A)$ depends on $(x^\perp,\theta)$,  where 
$x^\perp=x-(x\cdot\theta)\theta$.
Matrix $S(A)$ is called the non-abelian Radon transform of
$A(x,\theta)=A_1(x)\theta_1+A_2(x)\theta_2+A_0(x).$

In the abelian case $m=1$,  i.e.  when $c_0(x,\theta)$ and
$A_j(x),\ 0\leq j\leq 2$, are scalar functions we have an explicit
solution of (\ref{eq:1.1}):
\begin{equation}                             \label{eq:1.2}
c_0(x,\theta)=\exp\left(\int_{-\infty}^{x\cdot \theta}
A(x_\perp+s\theta,\theta)ds
\right).
\end{equation}
Therefore $S(A)(x_\perp,\theta)= \exp R(A)$,
where $R(A)(x_\perp,\theta)
=\int_{-\infty}^\infty A(x_\perp +s\theta,\theta)dt$
is the ordinary Radon transform of $A(x,\theta)$ (see [Na])
\qed

We shall 
call matrices $A^{(1)}(x,\theta)$ and $A^{(2)}(x,\theta)$,
$A^{(i)}(x,\theta)=\Sigma_{j=1}^2A_j^{(i)}(x,\theta)\theta_j+
A_0^{(i)}(x),\ i=1,2$  gauge equivalent if there exists 
a nonsingular $C^\infty$
matrix $g(x)$ such that $g(x)=I_m$ for $|x|\geq R$ and
\begin{eqnarray}                               \label{eq:1.3}
A_j^{(2)}=gA_j^{(1)}g^{-1}+\frac{\partial g}{\partial x_j} g^{-1},
\ \ j=1,2,
\\
A_0^{(2)}=gA_0^{(1)}g^{-1}.
\nonumber
\end{eqnarray}

If $c_0^{(1)}(x,\theta)$ satisfies (\ref{eq:1.1}) with $A(x,\theta)$
replaced by $A^{(1)}(x,\theta)$, then $c_0^{(2)}(x,\theta)=
g(x)c_0^{(1)}(x,\theta)$ satisfies 
\[
\theta\cdot \frac{\partial c_0^{(2)}}{\partial  x}=
A^{(2)}(x,\theta)c_0^{(2)}.
\]
Since $g(x)=I$ for $|x|>R$ we have that
$\lim_{s\rightarrow -\infty}c_0^{(2)}(x+s\theta,\theta)=I_m$ and
$\lim_{s\rightarrow +\infty}c_0^{(2)}(x+s\theta,\theta)
=S(A^{(1)})(x_\perp,\theta)$.
Therefore $A^{(1)}(x,\theta)$ and $A^{(2)}(x,\theta)$ 
have the same non-abelian Radon transform.
We shall prove an inverse statement:
\begin{theorem}                         \label{theo:1.1}
Suppose $A_j^{(1)}(x)$ and $A_j^{(2)}(x),\ 1\leq j\leq 2,$
are $C^\infty$ compactly supported matrices with the same
non-abelian Radon transform.  Then $A_j^{(1)}(x),\ 0\leq j\leq 2$ 
and $A_j^{(2)}(x),\ 0\leq j\leq 2,$  are gauge equivalent.
\end{theorem} 

In the case $A_0^{(1)}=A_0^{(2)}=0$ 
this result is contained in [E1]  but it was not explicitely stated 
there. In this paper we shall prove Theorem \ref{theo:1.1} in the case 
$A_0^{(j)}(x)\neq 0,\ j=1,2,$ 
and consider some extensions and applications.

The first works on the non-abelian Radon transform were done by
Wertgeim [We] and Sharifutdinov [Sh].  They proved the  uniqueness
of the inverse problem modulo gauge transformations assuming that
$A^{(1)}$ and $A^{(2)}$ 
are small.  A major work on this subject belongs to R.Novikov [N1].
He proved the global uniqueness modulo gauge transformations in
$\R^n,n\geq 3$.  In the case $n=2$ he also assumed that
$A(x,\theta)$ is small but he gave the reconstruction procedure 
using the Riemann-Hilbert problem.  Finch and Uhlmann [FU] proved
the uniqueness assuming that $A_0(x)=0,\ A_j(x),\ 1\leq j\leq 2,$
have compact supports and the curvature 
$\Omega_{12}=\frac{\partial A_1}{\partial x_2}
-\frac{\partial A_2}{\partial x_1}+[A_1,A_2]$ is small.
R.Novikov [N1] discovered examples of the non-uniqueness in
the non-abelian Radon transform : there are $A_1(x),A_2(x)$
with noncompact supports that are not gauge equivalent to the zero 
matrices and such that the corresponding non-abelian Radon transform
is $I_m$.  These examples 
appeared in the works of physicists [Wr], [V] on the theory
of solutions.  The non-uniqueness examples show that the global
uniqueness result given by Theorem \ref{theo:1.1} is not trivial.

Theorem \ref{theo:1.1} will be proved in \S 2.  The crucial part
of the proof is Lemma \ref{lma:2.1} proven in [ER3].  This lemma 
allows also
to extend to the non-abelian case the Novikov's formula for the
inversion of the attenuated Radon transform [N].  This will be done
in \S 3.  In \S 4 we apply the non-abelian Radon transform to the
inverse scattering problem for the Schr\"{o}dinger equation with
Yang-Mills potentials in two dimension.

\section{The proof of Theorem \ref{theo:1.1}.}
\label{section 2}
\init

Consider the equation:
\begin{equation}                             \label{eq:2.1}
\zeta_1\frac{\partial c}{\partial x_1} +
\zeta_2\frac{\partial c}{\partial x_2}=
(A_1(x)\zeta_1+A_2(x)\zeta_2 +A_0(x))c(x,t),
\end{equation}
where $A_j(x),\ 0\leq j\leq 2,$ are the same as in (\ref{eq:1.1}),
$\zeta_p\in \C,\ p=1,2,\ \zeta_1^2+\zeta_2^2 =1$.  Define
\[
\zeta_1(t)=\frac{1}{2}(t+ \frac{1}{t}),\ \ \ \ 
\zeta_2(t)=\frac{1}{2}(t- \frac{1}{t}), 
\]
where $t\in \C\setminus\{0\}$.
Denote $D^+=\{t:\ |t|>1\},\ \ D^-=\{t:\ |t|<1\}$.
The following lemma proven in [ER3] is the main part of the proof
of Theorem \ref{theo:1.1}:
\begin{lemma}                                \label{lma:2.1}
There exist solutions $c_\pm(x,t)$ of  
(\ref{eq:2.1}) with $\zeta_p=\zeta_p(t),\ p=1,2,$ having the following 
properties:

a) $c_+(x,t)$ and $c_-(x,t)$ are solutions of (\ref{eq:2.1})
for $(x,t)\in \overline{B_{2R}}\times\overline{D^+}$
and $(x,t)\in \overline{B_{2R}}\times\overline{D^-}$
respectively.

b) $c_+(x,t)\ (c_-(x,t))$ is smooth when
 $(x,t)\in \overline{B_{2R}}\times\overline{D^+}$ 
($ \overline{B_{2R}}\times\overline{D^-})$
and $\det c_+(x,t)\neq 0 \ (\det c_-(x,t)\neq 0)$ for all
$(x,t)\in \overline{B_{2R}}\times\overline{D^+}$ 
($ \overline{B_{2R}}\times\overline{D^-})$.

c) $c_+(x,t) (c_-(x,t))$ is analytic in $t$ when 
$(x,t)\in \overline{B_{2R}}\times\overline{D^+}$ 
($ \overline{B_{2R}}\times\overline{D^-})$.
Moreover $c_+(x,t)$ is analytic at $t=\infty$ with
$\det c_+(x,\infty)\neq 0$ and
$c_-(x,t)$ is analytic at $t=0$ with
$\det c_-(x,0)\neq 0$.
\end{lemma}

Note that
\begin{equation}                       \label{eq:2.2}
\zeta(t)\cdot\frac{\partial}{\partial x}=
\zeta_1(t)\frac{\partial}{\partial x_1}+
\zeta_2(t)\frac{\partial}{\partial x_2}=
t\frac{\partial}{\partial\bar{z}} + t^{-1}\frac{\partial}{\partial z},
\end{equation}
where $\frac{\partial}{\partial \bar{z}}
=\frac{1}{2}\left(\frac{\partial}{\partial x_1}
+i\frac{\partial}{\partial x_2}\right),\ 
\frac{\partial}{\partial z}
=\frac{1}{2}\left(\frac{\partial}{\partial x_1}
-i\frac{\partial}{\partial x_2}\right)$.

The operator $\zeta(t)\cdot \frac{\partial}{\partial x}$ is elliptic
operator when $|t|\neq 1$ and $\zeta(t)\cdot\frac{\partial}{\partial x}$
degenerates to $\theta(\varphi)\cdot\frac{\partial}{\partial x}$
when $|t|=1,\ t=e^{i\varphi},\ \theta(\varphi)=(\cos\varphi,-sin\varphi)$.
This makes the proof of Lemma \ref{lma:2.1} quite complicate (see [ER3])
\qed

Suppose $A^{(1)}(x,\theta)$ and $A^{(2)}(x,\theta)$ are such that 
the non-abelian transforms $S(A^{(1)})$ and $S(A^{92)})$ are equal.  Let 
$c_\pm^{(i)}(x,t)$ be the solutions of the equations
\begin{equation}                                \label{eq:2.3}
\zeta(t)\cdot \frac{\partial}{\partial x}c_\pm^{(i)}(x,t)=
A^{(i)}(x,\zeta(t))c_\pm^{(i)}(x,t),\ i=1,2,
\end{equation}
obtaines by the Lemma \ref{lma:2.1}.

Denote
\[
y_1=x\cdot\theta(\varphi),\ \ \ y_2=x\cdot\nu(\varphi),
\]
where $\nu=(\sin\varphi,\cos\varphi),$ and denote by
$c_\pm^{(i)}(y_1,y_2,\varphi)$ the matrix $c_\pm^{(i)}(x,e^{i\varphi})$ 
in $(y_1,y_2)$-coordinates.  Since  
(\ref{eq:2.3}) has the form 
\begin{equation}                              \label{eq:2.4}
\frac{\partial}{\partial y_1}c_\pm^{(i)}(y_1,y_2,\varphi)=
A^{(i)}(x(y,\varphi),\theta(\varphi))c_\pm^{(i)}(y_1,y_2,\varphi)
\end{equation}
in new coordinates when $t=e^{i\varphi}$ and since
$c_\pm^{(i)}(y_1,y_2,\varphi)(c_\pm^{(i)}(-\infty,y_2,\varphi))^{-1}$
is the solution of (\ref{eq:2.4}) equal to $I_m$ when 
$s\rightarrow\-\infty$ we get that
\begin{equation}                                 \label{eq:2.5}
S(A^{(i)})
=c_\pm^{(i)}(+\infty,y_2,\varphi)(c_\pm^{(i)}(-\infty,y_2,\varphi))^{-1},
\ \ i=1,2.
\end{equation}
Therefore $S(A^{(1)})=S(A^{(2)})$ implies that   
\begin{equation}                               \label{eq:2.6}                    
       c_\pm^{(1)}(+\infty,y_2,\varphi)
(c_\pm^{(1)}(-\infty,y_2,\varphi))^{-1}=
c_\pm^{(2)}(+\infty,y_2,\varphi)
(c_\pm^{(2)}(-\infty,y_2,\varphi))^{-1}.
\end{equation}
It follows from (\ref{eq:2.6}) that
\begin{equation}                                 \label{eq:2.7}
(c_+^{(2)}(+\infty,y_2,\varphi))^{-1}c_+^{(1)}(+\infty,y_2,\varphi)=
(c_+^{(2)}(-\infty,y_2,\varphi))^{-1}c_+^{(1)}(-\infty,y_2,\varphi),
\end{equation}
\begin{equation}                                 \label{eq:2.8}
(c_-^{(2)}(+\infty,y_2,\varphi))^{-1}c_-^{(1)}(-\infty,y_2,\varphi)=
(c_-^{(2)}(-\infty,y_2,\varphi))^{-1}c_-^{(1)}(-\infty,y_2,\varphi).
\end{equation}

Denote $Q_+(x,t)=(c_+^{(2)}(x,t))^{-1}c_+^{(1)}(x,t)$ when
$|x|>R,\ |t|\geq 1$.  Since $A^{(j)}(x,\zeta(t))=0$ for $|x|>R,\ j=1,2,$
we have that $\zeta(t)\cdot\frac{\partial}{\partial x}c_+^{(j)}(x,t)=0$
for $|x|>R,\ |t|\geq 1,\ j=1.2,$ and therefore
\begin{equation}                                \label{eq:2.9}
\zeta(t)\cdot\frac{\partial Q_+(x,t)}{\partial x}=0
\end{equation}
for $|x|>R,\ |t|\geq 1$.  We shall prove the following lemma:
\begin{lemma}                                     \label{lma:2.2}
Assume that (\ref{eq:2.7}) holds.  Then there exists a matrix
$Q_+(x,t)$, defined on $\overline{B_{2R}}\times\overline{D^+}$,
such that $Q_+(x,t)=(c_+^{(2)}(x,t))^{-1}c_+^{(1)}(x,t)$ for
$|x|>R,\ |t|\geq 1$ and has the following properties:

$a_1)$ $Q_+(x,t)\in C^\infty$  for 
$x\in \overline{B_{2R}},\ t\in \overline{D^+}$,

$b_1)$ $Q_+(x,t)$ is analytic in $t$ for $t\in D^+$ including $t=\infty$,

$c_1)$ $\det Q_+(x,t)\neq 0$ for 
$x\in \overline{B_{2R}},\ t\in \overline{D^+}$,
including $t=\infty$,

$d_1)$  Equation (\ref{eq:2.9}) holds for all $x\in \overline{B_{2R}},\ 
|t|\geq 1$.
\end{lemma}

Denote 
\begin{equation}                               \label{eq:2.10}
\Pi(t)f=\frac{1}{(2\pi)^2}
\int_{\R^2}\frac{\tilde{f}(\xi)e^{Ix\cdot\xi}d\xi_1d\xi_2}
{i(\zeta_1(t)\xi_1+\zeta_2(t)\xi_2)}
=\frac{1}{\pi}
\int_{\R^2}\frac{f(x_1',x_2')dx_1'dx_2'}
{t(z-z')+t^{-1}(\bar{z}-\bar{z}')}\ ,
\end{equation}
where $z=x_1+ix_2,z'=x_1'+ix_2'$.  Note that $\Pi(t)$ is the inverse
of $\zeta(t)\cdot\frac{\partial}{\partial x}$,  i.e. 
$\zeta(t)\cdot \frac{\partial}{\partial x}\Pi(t)f=f,
\ \forall f \in C_0^\infty(\R^2)$.
Denote
\begin{equation}                              \label{eq:2.11}
P_\pm^{(1)}(x,t)= c_\pm^{(1)}(x,t)-\Pi(t)A(x,\zeta(t))c_\pm^{(1)}(x,t),
\end{equation}
where 
$(x,t)\in \overline{B_{2R}}\times\overline{D^+}$
for $P_+^{(1)}$ and 
$(x,t)\in \overline{B_{2R}}\times\overline{D^-}$
for $P_-^{(1)}(x,t)$.  Applying $\zeta(t)\cdot\frac{\partial}{\partial x}$
to (\ref{eq:2.11}) we get
\[
\zeta(t)\cdot\frac{\partial}{\partial x}P_\pm^{(1)}(x,t)=0,\ \ \ j=1,2,
\]
where $(x,t)\in \overline{B_{2R}}\times\overline{D^+}$
for $P_+^{(1)}$ and 
$(x,t)\in \overline{B_{2R}}\times\overline{D^-}$
for $P_-^{(1)}$.

It was proven in [ER3](see the Remark in [ER3]) that 
$P_\pm^{(1)}(x,t)=g_\pm^{(1)}(\zeta^\perp\cdot x,t)$ where
$\zeta^\perp(t)=(-\zeta_2(t),\zeta_1(t)),\ g_+^{(1)}(z,t)$ is
analytic in $z,\ z\in \overline{B_{2R}}$
for $1\leq t\leq 2$ and $g_-^{(1)}(z,t)$ is
analytic in $z,\ z\in \overline{B_{2R}},\ \frac{1}{2}\leq|t|\leq 1.$
Note that $\zeta^\perp(e^{i\varphi})\cdot x=
\nu\cdot x=y_2$.  Therefore
$P_\pm^{(1)}(x,e^{i\varphi})=P_\pm^{(1)}(y_2,\varphi)$ is analytic
in $y_2\in\C$ for $|y_2|\leq 2R$.

Denote 
\[
\Pi_+(e^{i\varphi})=\lim_{r\rightarrow 1+0}\Pi(re^{i\varphi}),
\ \Pi_-(e^{i\varphi})=\lim_{r\rightarrow 1-0}\Pi(re^{i\varphi}).
\]
Here $r\rightarrow 1-0$ means $r\rightarrow 1$ and $r<1$, 
$r\rightarrow 1+0$
means that $r\rightarrow 1,\ r>1$.

It is easy to show
(see,for example, [ER1], formula (4.11) ) that
\begin{equation}                             \label{eq:2.12}
\Pi_+(e^{i\varphi})f=\int_{-\infty}^{y_1}(\Pi^+f)dy_1'
-\int_{y_1}^\infty(\Pi^-f)dy_1',
\end{equation}
\begin{equation}                             \label{eq:2.13}
\Pi_-(e^{i\varphi})f=\int_{-\infty}^{y_1}(\Pi^-f)dy_1'
-\int_{y_1}^\infty(\Pi^+f)dy_1',
\end{equation}
where 
\[
\Pi^\pm f=\frac{\mp i}{2\pi}\int_{-\infty}^\infty
\frac{f(y_1,y_2')dy_2'}{y_2-y_2'\mp i0},\ \ \ f\in C_0^\infty(\R^2).
\]
Passing in (\ref{eq:2.11}) to the limit when $r\rightarrow 1+0$ we get
for $j=1$:
\begin{equation}                            \label{eq:2.14}
P_\pm^{(1)}(y_2,\varphi)=c_\pm^{(1)}(y_1,y_2,e^{i\varphi})
-\Pi_\pm(e^{i\varphi})A^{(1)}(x(y,\varphi),\theta(\varphi))c_\pm^{(1)}.
\end{equation}
Taking the limits when $y_1\rightarrow\pm\infty$ and using
(\ref{eq:2.12}), (\ref{eq:2.13}) we get
\begin{eqnarray}                                   \label{eq:2.15}
c_\pm^{(1)}(+\infty,y_2,\varphi)=P_\pm^{(1)}(y_2,\varphi)
+\int_{-\infty}^\infty(\Pi^\pm(A^{(1)}c_\pm^{(1)})dy_1',
\\
c_\pm^{(1)}(-\infty,y_2,\varphi)=P_\pm^{(1)}(y_2,\varphi)
-\int_{-\infty}^\infty(\Pi^\mp(A^{(1)}c_\pm^{(1)})dy_1'.
\nonumber
\end{eqnarray}
Note that
$(c_\pm^{(2)}(x,t))^{-1}$ satisfies the equation
\[
\zeta(t)\cdot\frac{\partial}{\partial x}(c_\pm^{(2)})^{-1}=
-(c_\pm^{(2)})^{-1}A^{(2)}(x,\zeta(t)).
\]
Denote
\[
P_\pm^{(2)}(x,t)= (c_\pm^{(2)}))^{-1}
+\Pi(t)(c_\pm^{(2)})^{-1}A^{(2)}(x,\zeta(t)).
\]
Note that $\zeta(t)\cdot \frac{\partial}{\partial x}P_\pm^{(2)}=0$
for $|x|\leq 2R,\ t\in \overline{D^+}$ and $\overline{D^-}$,
respectively.  Since $(c_\pm^{(2)}(x,t))^{-1}$ has the same analytic 
properties as $c_\pm^{(2)}(x,t)$ the Remark in [ER3] applies to
$P_\pm^{(2)}(x,t)$.  In particular we have that
$P_\pm^{(2)}(x,e^{i\varphi})=P_\pm^{(2)}(y_2,\varphi)$ is analytic
for $y_2\in \C$ when $|y_2|\leq 2R.$

Taking the limits when $y_1\rightarrow\pm\infty$  we get
\begin{eqnarray}                                \label{eq:2.16}
\ \ \ \ \ \ \ \ \ 
(c_\pm^{(2)}(+\infty,y_2,\varphi))^{-1}=
P_\pm^{(2)}(y_2,\varphi)
-\int_{-\infty}^\infty\Pi^\pm(c_\pm^{(2)})^{-1}A^{(2)}dy_1',
\\
\ \ \ \ \ \ \ \ \
(c_\pm^{(2)}(-\infty,y_2,\varphi))^{-1}=
P_\pm^{(2)}(y_2,\varphi)
+\int_{-\infty}^\infty\Pi^\mp(c_\pm^{(2)})^{-1}A^{(2)}dy_1'.
\nonumber
\end{eqnarray}
Substituting (\ref{eq:2.15}) and (\ref{eq:2.16})
into (\ref{eq:2.7}) we get
\begin{eqnarray}                           \label{eq:2.17}
\ \ \ \ \ \ \ \ 
\\
(P_+^{(2)}(y_2,\varphi)
-\int_{-\infty}^\infty\Pi^+(c_+^{(2)})^{-1}A^{(2)}dy_1')
(P_+^{(1)}(y_2,\varphi)+
\int_{-\infty}^\infty(\Pi^+A^{(1)}c_+^{(1)})dy_1')
\nonumber
\\
\nonumber
=
(P_+^{(2)}(y_2,\varphi)
+\int_{-\infty}^\infty\Pi^-(c_+^{(2)})^{-1}A^{(2)}dy_1')
(P_+^{(1)}(y_2,\varphi)-
\int_{-\infty}^\infty(\Pi^-A^{(1)}c_+^{(1)})dy_1').
\end{eqnarray}

Since $\Pi^+f$ is analytic in $y_2$ for $\Im y_2<0$ and 
$\Pi^-f$ is analytic in $y_2$ for $\Im y_2>0$ we have
that the left hand side of (\ref{eq:2.17}) is analytic in $y_2$
for $\Im y_2< 0$ and $|y_2|\leq 2R$ and
the right hand side of (\ref{eq:2.17}) is analytic in $y_2$
for $\Im y_2> 0$ and $|y_2|\leq 2R$.  Denote
\begin{eqnarray}                              \label{eq:2.18}
Q_+(y_2,\varphi)=
(c_+^{(2)}(-\infty,y_2,\varphi))^{-1}c_+^{(1)}(-\infty,y_2,\varphi)
\\
\nonumber
= 
(c_+^{(2)}(+\infty,y_2,\varphi))^{-1}c_+^{(1)}(+\infty,y_2,\varphi).
\end{eqnarray}
Then $Q_+(y_2,\varphi)$ is analytic in $y_2$ in the disk $|y_2|<2R$.
In particular,  $Q_+(x,e^{i\varphi})$ is real analytic in $(x_1,x_2)$ 
for $x_1^2+x_2^2\leq 4R^2$ since $y_2=x\cdot\nu$.

Define 
\begin{equation}                              \label{eq:2.19}
\hat{Q}_+(x,t)=\frac{1}{2\pi i}
\int_{|t'|=1}\frac{Q_+(x,t')tdt'}{t'(t-t')}.
\end{equation}
Then $\hat{Q}_+(x,t)$ is analytic in $t$ for $|t|>1,\ |x|\leq 2R,\ 
\hat{Q}(x,t)$ is real analytic in $x$ for $|x|\leq 2R,\ |t|>1$.
Denote by $\hat{Q}_+(x,e^{i\varphi})$ the limit of $\hat{Q}_+(x,t)$
when $r\rightarrow 1+0,$  where $t=re^{i\varphi}$. We shall show
that $\hat{Q}_+(x,e^{i\varphi})=Q_+(x,e^{i\varphi})$.

When $|x|>R\ \ Q_+(x,t)=(c_+^{(2)}(x,t))^{-1}c_+^{(1)}(x,t)$ is
analytic for $|t|>1$ including $t=\infty$ and smooth for $|t|\geq 1$.
Therefore
\[
Q_+(x,t)=\frac{1}{2\pi i}
\int_{|t'|=1}\frac{Q_+(x,t')tdt'}{t'(t-t')}
\]
by the Cauchy formula.  Therefore
\[
Q_+(x,e^{i\varphi})=\hat{Q}_+(x,e^{i\varphi}) \mbox{for\ }\ \ |x|\geq R.
\]
Since $Q_+(x,e^{i\varphi})$ and $\hat{Q}_+(x,e^{i\varphi})$
are real analytic in $x$ for $|x|\leq 2R$ we get that 
$Q_+(x,e^{i\varphi})=\hat{Q}_+(x,e^{i\varphi})$ for all $|x|\leq 2R$.
Therefore $\hat{Q}_+(x,t)$ is tha analytic continuation of 
$Q_+(x,e^{i\varphi})$
in $t$ from $|t|=1$ to $D^+$.  We shall write from now on $Q(x,t)$ 
instead of $\hat{Q}_+(x,t)$.

We already proved that $Q_+(x,t)$ satisfies conditions $a_1)$ and $b_1)$.
Since $Q_+(x,t)$ is real analytic in $x$ for $|x|<2R$ and 
$Q_+(x,t)=(c_+^{(2)}(x,t))^{-1}c_+^{(1)}c_+^{(1)}(x,t)$ satisfies
(\ref{eq:2.9}) for $|x|>R,\ |t|\geq 1$ we get that (\ref{eq:2.9}) is 
satisfied for all $|x|\geq 2R$.  It remains to prove that 
$\det Q_+(x,t)\neq 0$ for $|x|\leq 2R,\ |t|\geq 1$.
Since $Q_+(x,e^{i\varphi})=Q_+(y_2,\varphi)$ is independent of
$y_1$
we have $\det Q_+(y_2,\varphi)\neq 0$ by taking $|y_1|>R.$  
When $|t|>1$ denote
$S(t)=\{z: z=x\cdot\zeta^\perp(t),|x|\leq 2R\}$.
Since (\ref{eq:2.9}) holds $Q_+(x,t)=h(x\cdot\zeta^\perp(t),t),$
where $h(z,t)$ is analytic in $z$ for $z\in S(t).$
We have
\[
Q_+(x,t)=(c_+^{(2)}(x,t))^{-1}c_+^{(1)}(x,t)
\]
 for $z\in \partial S(t)$.
Since $\det c_+^{(j)}(x,t)\neq 0$ for $|x|\leq 2R,\ j=1,2,$
we have that the increment of the argument of $\det h(z,t)$ on 
$\partial S(t)$ is zero.  Therefore  $\det h(z,t)$ has no zero for
$z\in S(t),\ |t|>1,$   i.e.  $c_1)$ holds.
\qed

Denote
\begin{equation}                                  \label{eq:2.20}
c_+^{(3)}(x,t)=c_+^{(2)}(x,t)Q_+(x,t).
\end{equation}
Then $c_+^{(3)}(x,t)$ satisfies (\ref{eq:2.3}) for $j=2$ and has 
the same properties as $c_+^{(2)}(x,t)$.  It follows from (\ref{eq:2.20})
and (\ref{eq:2.18}) that
\begin{equation}                               \label{eq:2.21}
c_+^{(3)}(\pm\infty,y_2,\varphi)=c_+^{(1)}(\pm\infty,y_2,\varphi).
\end{equation}
Analogously we can construct a matrix $Q_-(x,t)$ that extends 
$(c_-^{(2)}(x,t))^{(-1}c_-^{(1)}(x,t)$ from $|x|>R,\ |t|\leq 1$
to $(x,t)\in\overline{B_{2R}}\times \overline{D^-}$ and satisfies 
 $a_1),\ b_1),\ c_1),\ d_1)$  for 
$x\in \overline{B_{2R}},\ t\in \overline{D^-}$.  In particular,
\begin{equation}                            \label{eq:2.22}
Q_-(y_2,\varphi)=(c_-^{(2)}(-\infty,y_2,\varphi))^{-1}
c_-^{(1)}(-\infty,y_2,\varphi)
=
(c_-^{(2)}(+\infty,y_2,\varphi))^{-1}
c_-^{(2)}(+\infty,y_2,\varphi)
\end{equation}
Replace $c_-^{(2)}(x,t)$ by
\begin{equation}                            \label{eq:2.23}
c_-^{(3)}(x,t)=c_-^{(2)}Q_-(x,t).
\end{equation}
Then $c_-^{(3)}(x,t)$ satisfies (\ref{eq:2.3}) for 
$j=2,\ x\in \overline{B_{2R}},\ t\in \overline{D^-}$
and has the same properties as $c_-^{(2)}(x,t)$.  It follows from
(\ref{eq:2.22}), (\ref{eq:2.23})  that
\begin{equation}                               \label{eq:2.24}
c_-^{(3)}(\pm\infty,y_2,\varphi)=c_-^{(1)}(\pm\infty,y_2,\varphi).
\end{equation}
Denote
\[
g_1(x)=c_+^{(1)}(x,\infty),\ \ g_3(x)=c_+^{(3)}(x,\infty).
\]
We have $\det g_1(x)\neq 0,\ \det g_3(x)\neq 0$.  Since
$Q_+(x,t)=(c_+^{(2)}(x,t))^{(-1}c_+^{(1)}(x,t)$ for 
$|x|>R,\ |t|\geq 1$ we have
\begin{equation}                              \label{eq:2.25}
c_+^{(3)}(x,t)=c_+^{(1)}(x,t)
\end{equation}
for all $|t|\geq 1,\ |x|>R.$  In particular,
\begin{equation}                             \label{eq:2.26}
g_1(x)=g_3(x)\ \ \ \ \ \mbox{for\ } |x|> R.
\end{equation}
Replace $c_+^{(1)},\ c_-^{(1)}$ by 
\begin{equation}                             \label{eq:2.27}
\hat{c}_+^{(1)}(x,t)=g_1^{-1}(x)c_+^{(1)}(x,t),\ \
\hat{c}_-^{(1)}(x,t)=g_1^{-1}(x)c_-^{(1)}(x,t).
\end{equation}
Then $\hat{c}_\pm^{(1)}(x,t)$ satisfy
\begin{equation}                             \label{eq:2.28}
\zeta(t)\cdot\frac{\partial \hat{c}_\pm^{(1)}}{\partial x}=
\hat{A}^{(1)}(x,\zeta(t))\hat{c}_\pm^{(1)},
\end{equation}
where 
\begin{eqnarray}                               \label{eq:2.29}
\hat{A}^{(1)}(x)=g_1^{-1}(x)A_j^{(1)}(x)g_1(x)-g_1^{-1}(x)
\frac{\partial g_1}{\partial x},\ j=1,2,
\nonumber
\\
\hat{A}_0^{(1)}=g_1^{-1}(x)A_0{(1)}(x)g_1(x).
\end{eqnarray}
Analogously,  if we replace $c_\pm^{(3)}(x)$ by
\begin{equation}                             \label{eq:2.30}
\hat{c}_\pm^{(3)}(x,t)=g_3^{-1}(x)c_\pm^{(3)}(x,t),
\end{equation}
we get that $\hat{c}_\pm^{(3)}(x,t)$ satisfies 
\begin{equation}                             \label{eq:2.31}
\zeta(t)\cdot \frac{\partial \hat{c}_\pm^{(3)}}{\partial x}=
\hat{A}^{(3)}(x,\zeta(t))\hat{c}_\pm^{(3)},
\end{equation}
where 
\begin{eqnarray}                             \label{eq:2.32}
\hat{A}_j^{(3)}(x)=g_3^{-1}(x)A_j^{(2)}(x)g_3(x)-g_3^{-1}
\frac{\partial g_3}{\partial x},\ \ j=1,2,
\nonumber
\\
\hat{A}_0^{(3)}=g_3^{-1}A_0^{(2)}g_3.
\end{eqnarray}
Note that 
\begin{equation}                             \label{eq:2.33}
\hat{c}_+^{(1)}(x,\infty)=\hat{c}_+^{(3)}(x,\infty)=I_m.
\end{equation}
Denote
\begin{equation}                             \label{eq:2.34}
b_1(x,e^{i\varphi})=(\hat{c}_-^{(1)}(x,e^{i\varphi}))^{-1}
\hat{c}_+^{(1)}(x,e^{i\varphi}),
\end{equation}
\begin{equation}                             \label{eq:2.35}
b_3(x,e^{i\varphi})=(\hat{c}_-^{(3)}(x,e^{i\varphi}))^{-1}
\hat{c}_+^{(3)}(x,e^{i\varphi}).
\end{equation}
Since $\hat{c}_\pm^{(1)}$ satisfy the same equation (\ref{eq:2.28})
when $t=e^{i\varphi}$ we get 
\[
\frac{\partial}{\partial y_1}b_1(x,e^{i\varphi})=
-(\hat{c}_-^{(1)})^{-1}(\hat{A}^{(1)}\hat{c}_-^{(1)})(\hat{c}_-^{(1)})^{-1}
\hat{c}_+^{(1)}
+(\hat{c}_-^{(1)})^{-1}\hat{A}^{(1)}\hat{c}_+^{(1)}=0,
\]
i.e. $b_1(x,e^{i\varphi})$ is independent of $y_1$.  Therefore
\[
b_1(y_2,\varphi)=
(\hat{c}_-^{(1)}(-\infty,y_2,\varphi))^{-1}
\hat{c}_+^{(1)}(-\infty,y_2,\varphi).
\]
Analogously,  since $\hat{c}_\pm^{(3)}$ satisfy
the same equation (\ref{eq:2.31}) with $t=e^{i\varphi}$ we get that
\[
b_3(y_2,\varphi)=(\hat{c}_-^{(3)}(-\infty,y_2,\varphi))^{-1}
\hat{c}_+^{(3)}(-\infty,y_2,\varphi).
\]
It follows from (\ref{eq:2.21}), (\ref{eq:2.24}), (\ref{eq:2.27}),
(\ref{eq:2.30}) that 
\begin{equation}                               \label{eq:2.36}
b_1(x,e^{i\varphi})=b_3(x,e^{i\varphi}).
\end{equation}               
Therefore (\ref{eq:2.34}), (\ref{eq:2.35}), (\ref{eq:2.36}) imply
that
\begin{equation}                           \label{eq:2.37}
(\hat{c}_-^{(1)}(x,e^{i\varphi}))^{-1}
\hat{c}_+^{(1)}(x,e^{i\varphi})=
(\hat{c}_-^{(3)}(x,e^{i\varphi}))^{-1}
\hat{c}_+^{(3)}(x,e^{i\varphi}),
\end{equation}
for all $x\in B_{2R},\ \varphi\in[0,2\pi]$.

We can rewrite (\ref{eq:2.37}) as
\begin{equation}                             \label{eq:2.38}
\hat{c}_-^{(3)}(x,e^{i\varphi})
(\hat{c}_-^{(1)}(x,e^{i\varphi}))^{-1}
=\hat{c}_+^{(3)}(x,e^{i\varphi})
(\hat{c}_+^{(1)}(x,e^{i\varphi}))^{-1}.
\end{equation}
For each $x\in \overline{B_{2R}}$ the left hand side of (\ref{eq:2.38})
extends analytically to $D^-$ and the right hand side extends
analytically to $D^+$.  Therefore (\ref{eq:2.38}) defines an entire
matrix $d(t),\ t\in\C$.  It follows from (\ref{eq:2.33}) that $d(\infty)=I_m$.
Therefore by the Liouville theorem $d(t)=I_m$ for all $t$.  Therefore
\begin{equation}                             \label{eq:2.39}
\hat{c}_-^{(1)}(x,t)=\hat{c}_-^{(3)}(x,t)\ \ \mbox{for\ } |x|\leq 2R,\ |t|\leq 1.
\end{equation}
\begin{equation}                             \label{eq:2.40}
\hat{c}_+^{(1)}(x,t)=\hat{c}_+^{(3)}(x,t)\ \ \mbox{for\ } |x|\leq 2R,\ |t|\geq 1.
\end{equation}
Since 
\[
\hat{A}^{(j)}(x,\zeta(t))=
\left(\zeta(t)\cdot\frac{\partial \hat{c}_\pm^{(j)}}{\partial x}\right)
(\hat{c}_\pm^{(j)})^{-1},\ \ j=1,3.
\]
we get that
\begin{equation}                           \label{eq:2.41}
\hat{A}^{(1)}(x,\zeta(t))=\hat{A}^{(3)}(x,\zeta(t))
\end{equation}
for all $x\in \overline{B_{2R}}$ and $t\in \C\setminus \{0\}$.
Taking the limits when $t\rightarrow\infty$ and $t\rightarrow 0$ we get
\begin{eqnarray}                            \label{eq:2.42}
\hat{A}_1^{(1)}(x)+i\hat{A}_2^{(1)}(x)=\hat{A}_1^{(3)}(x)+i\hat{A}_2^{(3)}(x),
\\
\hat{A}_1^{(1)}(x)-i\hat{A}_2^{(1)}(x)=\hat{A}_1^{(3)}(x)-i\hat{A}_2^{(3)}(x).
\nonumber
\end{eqnarray} 
Therefore $\hat{A}_j^{(1)}(x)=\hat{A}_j^{(3)}(x),\ j=1,2.$
Then (\ref{eq:2.41}) implies that $\hat{A}_0^{(1)}(x)=\hat{A}_0^{(3)}(x).$
Therefore it follows from (\ref{eq:2.29}) and (\ref{eq:2.32}) that 
$A_j^{(1)}(x),\ 0\leq j\leq 2,$ and $A_j^{(2)}(x),\ 0\leq j\leq 2,$
are gauge equivalent with the gauge $g(x)=g_1^{-1}(x)g_2(x)$.
Note that $g(x)=I_m$ for $|x|>R$.
\qed

Note that (\ref{eq:2.34}), (\ref{eq:2.35})  are the Riemann-Hilbert 
problems on the circle $|t|=1$ for each $x\in\overline{B_{2R}}$.
The reduction of the inversion of the non-abelian Radon transform 
to the Riemann-Hilbert problem
was done first by R.Novikov in [N1] 
(see also [MZ] and [ER1], formulas (4.13), (4.14) )

\section{Attenuated non-abelian Radon transform.}
\label{section 3}
\init

Consider the following equation in $\R^2$:
\begin{equation}                              \label{eq:3.1}
\theta\cdot \frac{\partial u}{\partial x}-(A_1(x)\theta_1+A_2(x)\theta_2+
A_0(x))u(x,\theta)=f(x),
\end{equation}
where $A_j,\ 0\leq j\leq 2$ are smooth $m\times m$ matrices, 
$f(x),\ u(x,\theta)$ are $m$-vectors,  $\mbox{supp\ }A_j(x)\subset B_R,
\ 0\leq j\leq 2,\ \mbox{supp\ }f(x)\subset B_R$.

There is a unique solution
of (\ref{eq:3.1}) such that $u(x+s\theta,\theta)\rightarrow 0$ when 
$s\rightarrow -\infty$.  Let $c_0(x,\theta)$ be the same as in 
(\ref{eq:1.1}).  We look for $u(x,\theta)$ in the form
$u(x,\theta)=c_0(x,\theta)v(x,\theta)$.  Substituting in (\ref{eq:3.1})
we get
\begin{equation}                             \label{eq:3.2}
c_0(x,\theta)\theta\cdot\frac{\partial v(x,\theta)}{\partial x}= f(x).
\end{equation}
The unique solution of (\ref{eq:3.2}) such that 
$v(x+s\theta,\theta)\rightarrow 0$ when $s\rightarrow -\infty$ has the form
\[
v(x,\theta)=\int_{-\infty}^{x\cdot\theta}c_0^{-1}(x_\perp+\tau\theta,\theta)
f(x_\perp+\tau\theta)d\tau,
\]
where $x_\perp=x-(x\cdot\theta)\theta$.   Therefore
\begin{equation}                            \label{eq:3.3}
u(x,\theta)=c_0(x,\theta)\int_{-\infty}^{x\cdot\theta}
c_0^{-1}(x_\perp+\tau\theta,\theta)f(x_\perp+\tau\theta)d\tau.
\end{equation}
Take the limit
of $u(x+s\theta,\theta)$ when $s\rightarrow +\infty$.  We get
\[
\lim_{s\rightarrow +\infty}u(x+s\theta,\theta)=
\left(\lim_{s\rightarrow +\infty} c_0(x+s\theta,\theta)\right)R_A f,
\]
where
\begin{equation}                            \label{eq:3.4}
(R_Af)(x,\theta)=
\int_{-\infty}^\infty c_0^{-1}(x_\perp+\tau\theta,\theta)
f(x_\perp+\tau\theta)d\tau.
\end{equation}
The integral $R_Af$ is called
the attenuated Radon transform of $f(x)$.

We shall consider the inverse problem  
of recovering $f(x)$ knowing $R_Af$.  Matrices $A_j(x),\ 0\leq j\leq 2$,
are assumed to be known.

We shall repeat the reconstruction procedure of R.Novikov [N], however 
the result is new since it is based on the Lemma \ref{lma:2.1}.

Let $c_\pm(x,t)$ be the same matrices as in Lemma \ref{lma:2.1}.
Consider the following equation 
\begin{equation}                                \label{eq:3.5}
\zeta(t)\cdot\frac{\partial u_\pm(,t)}{\partial x}-A(x,\zeta(t))u_\pm(x,t)
=f(x),
\end{equation}
where $\zeta(t),\ A(x,\zeta(t))$ are the same as in (\ref{eq:2.3}),
$u_+(x,t)$ and $u_-(x,t)$ are defined on 
$\overline{B_{2R}}\times\overline{D^+}$ and 
$\overline{B_{2R}}\times\overline{D^-}$ respectively.  We look for
$u_\pm(x,t)$ in the form $u_\pm(x,t)=c_\pm(x,t)v_\pm(x,t)$. Then as
in (\ref{eq:3.2}) we get that $v_\pm(x,t)$ satisfy the following
equations:
\[
c_\pm(x,t)\zeta(t)\cdot\frac{\partial v_\pm(x,t)}{\partial x}=f(x).
\]
Therefore
\begin{equation}                                 \label{eq:3.6}
u_\pm(x,t)=c_\pm(x,t)\Pi(t)c_\pm^{-1}(x,t)f(x)
\end{equation}
are solutions of (\ref{eq:3.5}) in $\overline{B_{2R}}$,  where
$\Pi(t)$ is defined in (\ref{eq:2.10}).

Introduce coordinates $y_1=x\cdot\theta, \ y_2=x\cdot\nu$, 
where $\theta=(\cos\varphi,-\sin\varphi),\ \nu=(\sin\varphi,\cos\varphi),
\ t=re^{i\varphi}$.
Take the limit of $u_+(x,t)$ when $r\rightarrow 1+0$ and the limit
of $u_-(x,t)$ when $r\rightarrow 1-0$.

Denote as in \$ 2
\begin{equation}                                \label{eq:3.7}
u_\pm(y_1,y_2,\varphi)=u_\pm(x,e^{i\varphi}).
\end{equation}
Then we get
from (\ref{eq:3.6})
\begin{equation}                               \label{eq:3.8}
u_\pm(y_1,y_2,\varphi)=c_\pm(y_1,y_2,\varphi)\Pi_\pm(e^{i\varphi})
c_\pm^{-1}(y_1,y_2,\varphi)f(x(y,\varphi)),
\end{equation}
where
$\Pi_\pm(e^{i\varphi})$ have the form (\ref{eq:2.12}), (\ref{eq:2.13}).
Taking the limit when $y_1\rightarrow -\infty$ and
using (\ref{eq:2.12}),  (\ref{eq:2.13})
we get
\begin{equation}                              \label{eq:3.9}
u_+(-\infty,y_2,\varphi)=-c_+(-\infty,y_2,\varphi)
\int_{-\infty}^\infty(\Pi^-(c_+^{-1}f))(y_1',y_2)dy_1',
\end{equation}
\begin{equation}                              \label{eq:3.10}
u_-(-\infty,y_2,\varphi)=-c_-(-\infty,y_2,\varphi)
\int_{-\infty}^\infty(\Pi^+(c_-^{-1}f))(y_1',y_2)dy_1'.
\end{equation}
Note that $c_\pm(x,e^{i\varphi}),\ c_0(x,\theta)$ satisfy the
same homogeneous equation 
\[
\theta\cdot\frac{\partial c}{\partial x}=A(x,\theta)c.
\]
Therefore, by the uniqueness of the Cauchy problem
\begin{equation}                        \label{eq:3.11}
c_\pm(y_1,y_2,\varphi)=c_0(x,\theta(\varphi))c_\pm(-\infty,y_2,\varphi),
\end{equation}
since $c_0(x,\theta(\varphi))\rightarrow I_m$ when 
$y_1\rightarrow -\infty$.
Substituting (\ref{eq:3.11}) into (\ref{eq:3.9}) and (\ref{eq:3.10})
and taking into account that $\Pi^\pm$ commute with the integration
in $y_1$ we get:
\begin{equation}                              \label{eq:3.12}
u_-(-\infty,y_2,\varphi)=
-c_-(-\infty,y_2,\varphi)\Pi^+
\int_{-\infty}^\infty c_-^{-1}(-\infty,y_2,\varphi)(c_0^{-1}f)(y_1,y_2)dy_1,
\end{equation}

\begin{equation}                              \label{eq:3.13}
u_+(-\infty,y_2,\varphi)=
-c_+(-\infty,y_2,\varphi)\Pi^-
\int_{-\infty}^\infty c_+^{-1}(-\infty,y_2,\varphi)(c_0^{-1}f)(y_1,y_2)dy_1.
\end{equation}
Note that
\[
(R_Af)(y_2,\varphi)=\int_{-\infty}^\infty(c_0^{-1}f)(y_1,y_2)dy_1
\]
is known.
Therefore
\begin{equation}                              \label{eq:3.14}
u_\pm(-\infty,y_2,\varphi)=
-c_\pm(-\infty,y_2,\varphi)\Pi^\mp
 c_\pm^{-1}(-\infty,y_2,\varphi)(R_Af)(y_2,\varphi)
\end{equation}
are known too.

Since $u_+(x,e^{i\varphi})-u_-(x,e^{i\varphi})$ satisfies homogeneous 
equation $\theta\cdot \frac{\partial u}{\partial x}-A(x,\theta)u=0$
we have analogously to (\ref{eq:3.11}) that
\begin{equation}                             \label{eq:3.15}
u_+(x,e^{i\varphi})-u_-(x,e^{i\varphi})=c_0(x,\theta)
(u_+(-\infty,y_2,\varphi)-u_-(-\infty,y_2,\varphi)).
\end{equation}
Since $c_0(x,\theta)$ is known we can recover 
$u_+(x,e^{i\varphi})-u_-(x,e^{i\varphi})$.  Consider the integral
\begin{equation}                            \label{eq:3.16}
\frac{1}{2\pi i}\int_{|t|=1}(u_+(x,t)-u_-(x,t))dt.
\end{equation}
It follows fom (\ref{eq:3.6}) that for each $x\in\overline{B_{2R}}\ 
u_-(x,t)$ is analytic when $|t|<1$ and $u_-(x,t)$ is continuous when 
$|t|\leq 1$.
Therefore
\begin{equation}                             \label{eq:3.17}
\int_{|t|=1}u_-(x,t)dt=0.
\end{equation}
It follows also from (\ref{eq:3.6}) that $u_+(x,t)$ is analytic 
when $|t|>1$.  Note that
when $h(x)$ has a compact support and $t\rightarrow\infty$ we get
\begin{equation}                              \label{eq:3.18}
\Pi(t)h=\frac{1}{\pi}\int\int_{\R^2}
\frac{h(y_1,y_2)dy_1dy_2}{t(z-w)+\frac{1}{t}(\bar{z}-\bar{w})}=
\frac{1}{t}Sh+O\left(\frac{1}{t^2}\right),
\end{equation}
where $z=x_1+ix_2,\ w=y_1+iy_2$,
\begin{equation}                               \label{eq:3.19}
Sh=
\frac{1}{\pi}\int\int_{\R^2}
\frac{h(y_1,y_2)dy_1dy_2}{z-w}.
\end{equation}
Taking the limit in (\ref{eq:3.6}) when $t\rightarrow\infty$ we get
\[
u_+(x,t)=c_+(x,\infty)\frac{1}{t}S(c_+^{-1}(x,\infty)f(x)+
O\left(\frac{1}{t^2}\right).
\]
Therefore computing the residue at $t=\infty$ we get
\begin{equation}                            \label{eq:3.20}
\frac{1}{2\pi i}\int_{|t|=1|}u_+(x,t)dt
=c_+(x,\infty)S(c_+^{-1}(x,\infty)f(x)).
\end{equation}
Note that $S$ is the inverse to the Cauchy-Riemann operator 
$\frac{\partial}{\partial \bar{z}}=
\frac{1}{2}\left(\frac{\partial}{\partial x_1}
+i\frac{\partial}{\partial x_2}\right),$
i.e. $\frac{\partial}{\partial\bar{z}}Sh=h,
\ \ \forall\in C_0^\infty(\R^2)$.
Therefore multiplying (\ref{eq:3.20}) by $c_+^{-1}(x,\infty)$ 
from the left and then applying operator 
$\frac{\partial}{\partial\bar{z}}$
we have
\begin{equation}                                      \label{eq:3.21}
f(x)=c_+(x,\infty)\frac{\partial}{\partial\bar{z}}
\left(c_+^{-1}(x,\infty)\frac{1}{2\pi i}\int_{|t|=1}
(u_+(x,t)-u_-(x,t))dt\right).
\end{equation}
Therefore (\ref{eq:3.21}), (\ref{eq:3.15}), (\ref{eq:3.14}) give
the inversion formula for the attenuated Radon transform.
Note that $c_+(x,\infty)$ satisfies the equation 
\begin{equation}                                  \label{eq:3.22}
\left(\frac{\partial}{\partial x_1}+i\frac{\partial}{\partial x_2}\right)
c_+(x,\infty)=
(A_1(x)+iA_2(x))c_+(x,\infty)
\end{equation}
and $u_\pm(-\infty,y_2,\varphi)$ are given by (\ref{eq:3.14}).

Although matrices $c_\pm(x,t)$ are not unique any choice of 
$c_\pm(x,t)$ satisfying conditions of Lemma \ref{lma:2.1}
leads to a formula (\ref{eq:3.21}).

\section{Inverse scattering problem for the Schr\"{o}dinger equation with
exponentially decreasing Yang-Mills potentials.}
\label{section 4}
\init

Consider the following
equation in $\R^2$:
\begin{equation}                                    \label{eq:4.1}
\left[\sum_{j=1}^2\left(-i\frac{\partial}{\partial x_j}+A_j(x)\right)^2
+V(x)-k^2\right]u=0,
\end{equation}
where 
\begin{equation}                              \label{eq:4.2}
\left|\frac{\partial^p A_j(x)}{\partial x^p}\right|
\leq C_pe^{-\delta|x|},\ \ j=1,2,
\ \ \left|\frac{\partial^p V}{\partial x^p}\right|
\leq C_pe^{-\delta|x|},
\end{equation}
$\forall|p|\geq 0,\ \delta>0,\ A_j(x),j=1,2,\ V(x)$ and $u(x)$ are
$m\times m$ matrices.  Let $u(x)$ be a distorted plane wave in $\R^2$
having the following asymptotics:
\[
u=e^{ik\omega\cdot x}I_m+
\frac{a(\theta,\omega,k)e^{ik|x|}}{|x|^{\frac{1}{2}}}+
O\left(\frac{1}{|x|^{\frac{3}{2}}}\right),
\]
where $|x|\rightarrow\infty,\ \theta=\frac{x}{|x|},\ |\omega|=1,$ 
matrix $a(\theta,\omega,k)$ is the scattering amplitude.
The inverse scattering problem consists in the recovery of
$A_j(x), j=1,2,\ V(x)$ modulo gauge transformation,  knowing
the scattering amplitude $a(\theta,\omega,k)$.  Here the gauge
equivalence means (\ref{eq:1.3}) with $A_0(x)$ replaced by $V(x)$ 
and with $g(x)\in C^\infty(\R^2),\ \det g(x) \neq 0$ and 
$\lim_{|x|\rightarrow\infty}g(x)=I_m,\ \frac{\partial g}{\partial x}=
O(e^{-\delta|x|})$ when $|x|\rightarrow\infty.$

We assume that $A_j(x),\ j=1,2,$ satisfy the following 

\textbf{Condition (A)}\emph{ :There exist matrices $c_+(x,t)$ and 
$c_-(x,t)$
satisfying (\ref{eq:2.1}) in $\R^2\times\overline{D^+}$ and
$\R^2\times\overline{D^-}$  respectively with $A_0=0,\ A_1,A_2$ replaced 
by $-iA_1,-iA_2$, such that $\lim_{|x|\rightarrow\infty} c_\pm(x,t)=I_m$
and conditions a), b), c), d) of Lemma \ref{lma:2.1} are satisfied with
$\overline{B_{2R}}$ replaced by $\R^2$.}

Note that $c_\pm(x,t)$ satisfy the following equation in $\R^2$:
\begin{equation}                             \label{eq:4.3}
c_\pm(x,t)-i\Pi(t)A(x,\zeta(t))c_\pm=I_m.
\end{equation}

The Condition (A) allows to extend Theorem 2.2 of [E1] to the case
of potentials having noncompact supports.

\begin{theorem}                             \label{theo:4.1}
Suppose Condition (A) is satisfied.  Then knowing the scattering 
amplitude for all $k\in (k_0-\e,k_0+\e)$ and all 
$(\omega,\theta)\in S^1\times S^1$ 
we can recover $A_j(x), j=1,2,V(x)$ modulo a gauge transformation.
\end{theorem}
\underline{Proof}:
Since $A_j(x),j=1,2,V(x)$ are exponentially decreasing, the scattering 
amplitude $a(\theta,\omega,k)$ is real analytic in $k$.
Therefore $a(\theta,\omega,k)$ is known for all $k$ except possibly
a discrete set of $\{k_p\}_{p=0}^\infty$.

Let $c_\pm(x,t)$ be the matrices satisfying the Condition (A).
Then $c_\pm(x,t)$ are the solutions of 
\begin{equation}                           \label{eq:4.4}
\zeta(t)\cdot\frac{\partial c_\pm(x,t)}{\partial x}
=-iA(x)\cdot\zeta(t)c_\pm(x,t),
\end{equation}
where 
$c_+(x,t)\ \ (c_-(x,t))$ is defined in $\R^2\times \overline{D^+}\ \ 
(\R^2\times \overline{D^+}$),
$\zeta(t)=\frac{1}{2}(r+\frac{1}{r})\mu+i\frac{1}{2}(r-\frac{1}{r})\nu,
\ t=re^{i\varphi}, \ \mu=(\cos\varphi,-\sin\varphi),
\ \nu=(\sin\varphi,\cos\varphi).$

We shall define $c_+(x,\eta'+i\tau\nu)$ for all $x\in\R^2,\ \tau>0,
\eta'\in\R^2,\eta'\cdot\nu=0$ as follows:

$a_2)$ When $\eta'\cdot\mu>0$ and $(\eta'+i\tau\nu)^2=
|\eta'|^2-\tau^2>1$ we define 
$c_+(x,\eta'+i\tau)=c_+(x,t_1),$  where $|t_1|>1$ and
$\zeta(t_1)=\frac{\eta'+i\tau\nu}{(|\eta'|^2-\tau^2)^{\frac{1}{2}}}$.

$b_2)$ When $\eta'\cdot\mu<0,\ |\eta'|^2-\tau^2>1$, we  define
$c_+(x,\eta'+i\tau\nu)=c_-(x,t_2)$ where $|t_2|<1$ and
 $-\zeta(t_2)=\frac{\eta'+i\tau\nu}{(|\eta'|^2-\tau^2)^{\frac{1}{2}}}$.

$b_3)$  When $|\eta'|^2-\tau^2\leq 1,\ x\in \R^2$,  we define 
$c_+(x,\eta'+i\tau\nu)$ arbitrary requiring only that 
$c_+(x,\eta'+i\tau\nu)$ is smooth for all 
$(x,\eta',\tau)\in\R^1\times\R^1\times\R_+^1,\ 
c_+(x,\eta'+i\tau\nu)$
is homogeneous in $\eta'+i\tau\nu$ of degree zero,
$\det c_+(x,\eta'+i\tau\nu)\neq 0$ for all $(x,\eta'+i\tau\nu)$ and
$c_+(x,\eta'+i\tau\nu)\rightarrow I_m$ when $|x|\rightarrow \infty$.
Note that such $c_+(x,\eta'+i\tau\nu)$ exists since there is no
topological obstruction to the extention of $c_+(x,\eta'+i\tau\nu)$ from
$|\eta'|-\tau^2\geq 1,\ x\in\R^2$ to $|\eta'|-\tau\leq 1,\ x\in \R^2$,
where $c_+(x,\eta'+i\tau\nu)$ satisfies $c_+(\infty,\eta'+i\tau)=I_m,
\ \det c_+(x,\eta'+i\tau)\neq 0$.

Denote by
$c_+(x,D'+\xi+i\tau\nu)$ the pseudodifferential operator with the 
symbol $c_+(x,\eta'+\xi+i\tau\nu),$ where $\tau>0,\ \eta'\cdot\nu=0,
\ \xi=(k^2+\tau^2)^{\frac{1}{2}}\mu$.

Consider differential equation
\begin{equation}                                 \label{eq:4.5}
\left[\left(-\frac{\partial}{\partial x}
+\xi+i\tau\nu+A(x)\right)^2+V(x)-k^2\right]v=f.
\end{equation}
As in [ER], [ER1], [ER2] we are looking for
the solution of (\ref{eq:4.5}) in the form
\begin{equation}                                 \label{eq:4.6}
v=c_+(x,D'+\xi+i\tau\nu)Eg,
\end{equation}
where
\[
Eg=\frac{1}{(2\pi)^2}\int_{\R^2}
\frac{e^{ix\cdot\eta}\tilde{g}(\eta)d\eta}{(\eta+\xi+i\tau\nu)^2-k^2}.
\]
Substituting (\ref{eq:4.6}) into (\ref{eq:4.5}) we get
\begin{equation}                                 \label{eq:4.7}
c(x,D'+\xi+i\tau\nu)g+T_1g+T_2g=f,
\end{equation}
where
\[
T_1g=\frac{1}{(2\pi)^2}\int_{\R^2}
\frac{2(\eta+\xi+i\tau\nu)\cdot\left(A(x)c_+
-i\frac{\partial c_+}{\partial x}\right)
\tilde{g}(\eta)e^{ix\cdot\eta}d\eta}
{(\eta+\xi+i\tau\nu)^2-k^2},
\]
\[
T_2g=\frac{1}{(2\pi)^2}\int_{\R^2}
\frac{\left[\left(-\frac{\partial}{\partial x}+A(x)\right)^2+V(x)\right]
c_+\tilde{g}(\eta)e^{ix\cdot\eta}d\eta}
{(\eta+\xi+i\tau\nu)^2-k^2}.
\]

Note that $c_+$ is an elliptic pseudodifferenial operator.  In order
to solve (\ref{eq:4.7}) it is enough to show that the norms of $T_1$
and $T_2$ tends to zero when $\tau\rightarrow \infty,\ k\rightarrow\infty,
\ \frac{\tau}{k}\rightarrow 0$ (c.f. [ER], [ER1], [E]).
The estimates of $T_2$ is the same as in [ER], [ER1], [E].
Since 
$\eta'+\xi+i\tau\nu)\cdot\left(A(x)c_+
-i\frac{\partial c_+}{\partial x}\right)=0$
for $|\eta'+\xi|^2>\tau^2+1$.
We get
\[
T_1g=\frac{1}{(2\pi)^2}\int_{\R^2}
\frac{2(\eta\cdot\nu)
\left(A(x)c_+-i\frac{\partial c_+}{\partial x}\right)
\tilde{g}(\eta)e^{ix\cdot\eta}d\eta}
{(\eta+\xi+i\tau\nu)^2-k^2}
\]
\[
+\frac{1}{(2\pi)^2}\int_{\R^2}
\frac{\chi_-\left(\frac{|\eta'+\xi|^2}{\tau^2+1}\right)2(\eta'+\xi+i\tau\nu)
\left(A(x)c_+-i\frac{\partial c_+}{\partial x}\right)
\tilde{g}(\eta)e^{ix\cdot\eta}d\eta}
{(\eta+\xi+i\tau\nu)^2-k^2}
\]
\[
\stackrel{\mbox{def}}{=} T_{11}g+T_{12}g.
\]
Here $\eta'=\eta-(\eta\cdot\nu)\nu,\ \eta'\cdot \nu=0,
\ \chi_-(s)\in C^\infty(\R^1),\ \chi_-(s)=1$ for $s<1,\ \chi_-(s)=0$  
for $s\geq 2,\ \tau$ is large.  We have
\[
(\eta+\xi+i\tau\nu)^2-k^2=
(\eta +\xi)^2+2i\tau(\eta\cdot\nu)-\tau^2-k^2.
\]
Therefore
\[
\left|\Im \left[(\eta+\xi+i\tau\nu)^2-k^2\right]\right|\geq 
2\tau|(\eta\cdot\nu)|.
\]
This inequality implies that the symbol of $T_{11}$ is 
$O\left(\frac{1}{\tau}\right)$.
Also we have 
\[
\left|(\eta+\xi+i\tau\nu)^2-k^2\right|\geq
\frac{1}{2}\left|(\eta'+\xi)^2+\eta_\nu^2-\tau^2-k^2\right|
+\frac{1}{2}\tau|\eta_\nu|.
\]
When
$|\eta_\nu|\geq \tau$ we have $\left|(\eta+\xi+i\tau\nu)^2-k^2\right|\geq
\frac{1}{2}\tau^2$ and when $|\eta_\nu|\leq \tau$ and 
$|\eta'+\xi|\leq C\tau$ we have 
$\left|(\eta+\xi+i\tau\nu)^2-k^2\right|\geq
\frac{1}{2}k^2\geq\tau^2$ since $\frac{\tau}{k}\rightarrow 0$ when 
$\tau\rightarrow\infty$.  Therefore the norm of $T_{12}$ is 
$O\left(\frac{1}{\tau}\right)$.
The continuation of the proof of Theorem \ref{theo:4.1}
is the same as in [ER], [ER1] and [ER2].  In particular, we get
that the scattering amplitude $a(\theta,\omega,k)$  determines 
the following integral (see [ER2],  formula (2.29), 
or [ER1],  formula (4.8)):
\begin{equation}                             \label{eq:4.8}
I_+(y_2,\varphi)=-i\int_{_\infty}^\infty c_+^{-1}(y_1,y_2,\varphi)
(A(x)\cdot\mu(\varphi))dy_1,
\end{equation}
where $\theta(\varphi)=\mu(\varphi),\ y_1=\mu\cdot\ x,\ y_2=\nu\cdot x,\ 
c_+(y_1,y_2,\varphi)$ is $c_+(x,e^{i\varphi})$ in $(y_1,y_2)$-coordinates.

Define
matrix $c_-(x,\eta'-i\tau\nu)$, where $\tau>0$,  analogously to
$c_+(x,\eta'+i\tau\nu)$:

$a_3)$ $c_-(x,\eta'-i\tau\nu)=c_-(x,t_1')$,  where 
$\frac{\eta'-i\tau\nu}{|\eta'|^2-\tau^2)^{\frac{1}{2}}}=
\zeta(t_1'),\ |t_1'|<1$,  assuming that 
$\eta'\cdot\mu>0,\ |\eta'|^2-\tau^2>1$. 

$b_3)$ $c_-(x,\eta'-i\tau\nu)=c_+(x,t_2')$,  where 
$\frac{\eta'-i\tau\nu}{|\eta'|^2-\tau^2)^{\frac{1}{2}}}=
-\zeta(t_2'),\ |t_2'|<1$, 
assuming that 
$\eta'\cdot\mu<0,\ |\eta'|^2-\tau^2>1$. 

$c_3)$ $c_-(x,\eta'-i\tau\nu)$ has the same properties as
$c_+(x,\eta'+i\tau\nu)$ when $|\eta'|^2-\tau^2\leq 1$.

Repeating the proof with $c_-(x,\eta'+\xi-i\tau\nu)$ instead of
$c_+(x,\eta'+\xi +i\tau\nu)$ we get, taking the limit
when $k\rightarrow\infty,\ \tau\rightarrow \infty,
\ \frac{\tau}{k}\rightarrow 0$,
that the scattering amplitude determines the integral
\begin{equation}                               \label{eq:4.9}
I_-(y_2,\varphi)=
-i\int_{-\infty}^\infty c_-^{-1}(y_1,y_2,\varphi)(A\cdot\mu(\varphi))dy_1.
\end{equation}

Note that
$c_\pm^{-1}(y_1,y_2,\varphi)$ satisfy the equation
\[
-\frac{\partial}{\partial y_1}c_\pm^{-1}=-ic_\pm^{-1}(A\cdot\mu).
\]
Therefore
\begin{equation}                             \label{eq:4.10}
I_+(y_2,\varphi)=-\int_{-\infty}^\infty
\frac{\partial}{\partial y_1}c_+^{-1}dy_1=
c_+^{-1}(-\infty,y_2,\varphi)-c_+^{-1}(+\infty,y_2,\varphi).
\end{equation}
As in (\ref{eq:2.16}) we have, with $P_\pm^{(2)}$ replaced by $I_m$, that
\begin{equation}                             \label{eq:4.11}
c_+^{-1}(\pm\infty,y_2,\varphi)-I_m=
\mp i\int_{-\infty}^\infty \Pi^\pm(c_+)^{-1}(A\cdot\mu)dy_1.
\end{equation}
It follows from (\ref{eq:4.11}) that
\begin{equation}                             \label{eq:4.12}
c_+^{-1}(\pm\infty,y_2,\varphi)-I_m=\pm\Pi^\pm I_+.
\end{equation}
Therefore we can recover $c_+^{-1}(\pm\infty,y_2,\varphi)$.
Note that the recovery of $c_+^{-1}(\pm\infty,y_2,\varphi)$ from
(\ref{eq:4.8}) is the same as the computations (4.8)-(4.14) in [ER1].

Analogously starting from $I_-=\int_{-\infty}^\infty
-ic_-^{-1}(y_1,y_2,\varphi)(A\cdot\mu)dy_1$
we can recover $c_-^{-1}(\pm\infty,y_2,\varphi).$  Denote
\begin{equation}                                \label{eq:4.13}
b(x,e^{i\varphi})=\left(c_-(x,e^{i\varphi})\right)^{-1}c_+(x,e^{i\varphi}).
\end{equation}
Since $c_-$ and $c_+$ satisfy the same differential equation
\[
\frac{\partial}{\partial y_1}c=-i(A\cdot \mu)c,
\]
we get that $\frac{\partial}{\partial y_1} b(x,x^{i\varphi})=0$, i.e.
$b(x,e^{i\varphi})$ is independent of $y_1:\ b(x,e^{i\varphi})
=b(y_2,\varphi)$. 
Since we recovered $c_\pm(-\infty,y_2,\varphi)$ we know 
$b(y_2,\varphi)$:
\begin{equation}                             \label{eq:4.14}
b_2(y_2,\varphi)=\left(c_-(-\infty,y_2,\varphi)\right)^{-1}
c_+(-\infty,y_2,\varphi).
\end{equation}
We consider (\ref{eq:4.14}) for each $x\in \R^2$ as a Riemann-Hilbert
problem on the circle $|t|=1$ where $b(x,e^{i\varphi})$ is known and
$c_\pm(x,e^{i\varphi})$ are the unknowns.  If $c_\pm^{(1)}(x,t)$ is
another solution of the Riemann-Hilbert problem (\ref{eq:4.14}) then
we have that $\det c_+^{(1)}(x,t)\neq 0\ (\det c_-^{(1)}(x,t)\neq 0)$
for $(x,t)\in\R^2\times\overline{D^+}\ (\R^2\times\overline{D^-})$
respectively,  $\det c_+^{(1)}(x,\infty)\neq 0,
\ \det c_\pm^{(1)}(\infty,t)=I_m$.  Since
\begin{equation}                          \label{eq:4.15}
\left(c_-^{(1)}(x,e^{i\varphi})\right)^{-1}c_+^{(1)}(x,e^{i\varphi})
=\left(c_-(x,e^{i\varphi})\right)^{-1}c_+(x,e^{i\varphi}),
\end{equation}
we get by the Liouville theorem that 
\[
c_\pm^{(1)}(x,t)c_\pm^{-1}(x,t)=g(x),
\]
where $\det g(x)\neq 0,\ x\in \R^2,\ \det g(\infty)=I_m$.
Therefore $c_\pm^{(1)}(x,t)=g(x)c_\pm(x,t)$ satisfy the equation:
\[
\zeta(t)\cdot\frac{\partial c_\pm^{(1)}}{\partial}=
-i(A'(x)\cdot\zeta(t))c_\pm^{(1)}(x,t),
\]
where $A'(x)=(A_1'(x),A_2'(x))$ is gauge equivalent to
$A(x)=(A_1(x),A_2(x))$.
\qed

Now we shall recover $V(x)$ assuming that we already know $A(x)$
and $c_\pm(x,t)$.  As in [E] (see [E], formula (6.30), see also [ER],
formula (79),  and [ER2], formula (3.43) )  we get that the scattering
amplitude allows to recover
\begin{equation}                                       \label{eq:4.16}
J_+=\int_{-\infty}^\infty c_+^{-1}(y_1,y_2,\varphi)V(x(y,\varphi))
c_+(y_1,y_2,\varphi)dy_1
\end{equation}
and
\begin{equation}                                       \label{eq:4.17}
J_-=\int_{-\infty}^\infty c_-^{-1}(y_1,y_2,\varphi)V(x(y,\varphi))
c_-(y_1,y_2,\varphi)dy_1.
\end{equation}
Denote
\begin{equation}                                       \label{eq:4.18}
B_\pm(x,t)=c_\pm(x,t)\left(\Pi(t)c_\pm^{-1}(x,t)V(x)
c_\pm(x,t)\right)c_\pm^{-1},
\end{equation}
where $B_+(x,t)$ is defined on $\R^2\times \overline{D^+}$ and $B_-(x,t)$
is defined on $\R^2\times\overline{D^-}$.  We have
\begin{eqnarray}                                       \label{eq:4.19}
\zeta(t)\cdot\frac{\partial B_\pm}{\partial x}=
\left(\zeta(t)\cdot\frac{\partial c_\pm}{\partial x}\right)
\left(\Pi(t)(c_\pm^{-1}Vc_\pm)\right)c_\pm^{-1}
+c_\pm(c_\pm^{-1}Vc_\pm)c_\pm^{-1}
\\
+c_\pm\left(\Pi(t)(c_\pm^{-1}Vc_\pm)\right)\zeta(t)\cdot
\frac{\partial c_\pm^{-1}}{\partial x}
=-iA(x)\cdot\zeta(t)B_\pm+iB_\pm A(x)\cdot\zeta(t)+V(x).
\nonumber
\end{eqnarray}
Taking the limits
when $r\rightarrow 1+0$ and $r\rightarrow 1-0,\ t=re^{i\varphi}$
we get
\begin{equation}                                      \label{eq:4.20}
\mu(\varphi)\cdot\frac{\partial B_\pm(x,e^{i\varphi})}{\partial x}
=-iA(x)\cdot\mu(\varphi)B_\pm(x,e^{i\varphi})
+iB_\pm(x,e^{i\varphi})A(x)\cdot\mu(\varphi)+V(x).
\end{equation}
Introducing $(y_1,y_2)$ coordinates as before and using (\ref{eq:2.12})
and (\ref{eq:2.13}) we get
\[
B_\pm(-\infty,y_2,\varphi)=-c_\pm(-\infty,y_2,\varphi)
\int_{-\infty}^\infty\left(\Pi^\mp(c_\pm^{-1}Vc_-)\right)dy_1
c_\pm^{-1}(-\infty,y_2,\varphi).
\]
Using (\ref{eq:4.16}) and (\ref{eq:4.17}) we have
\begin{equation}                                   \label{eq:4.21}
B_+(-\infty,y_2,\varphi)
=-c_+(-\infty,y_2,\varphi)(\Pi^-J_+(y_2,\varphi))c_+^{-1}(-\infty,y_2,\varphi),
\end{equation}

\begin{equation}                                   \label{eq:4.22}
B_-(-\infty,y_2,\varphi)
=-c_-(-\infty,y_2,\varphi)(\Pi^+J_-(y_2,\varphi))c_-^{-1}(-\infty,y_2,\varphi),
\end{equation}
i.e. we can recover 
$B_\pm(-\infty,y_2,\varphi)$ from the scattering amplitude.  Consider
\begin{equation}                                    \label{eq:4.23}
B(x,e^{i\varphi})=B_+(x,e^{i\varphi})-B_-(x,e^{i\varphi}).
\end{equation}
Since $B_\pm(x,e^{i\varphi})$ satisfy the same equation (\ref{eq:4.20})
we get
\begin{equation}                                     \label{eq:4.24}
\mu(\varphi)\cdot\frac{\partial B}{\partial x}=
-iA(x)\cdot\mu(\varphi)B(x,e^{i\varphi})
+iB(x,e^{i\varphi})A(x)\cdot \mu(\varphi).
\end{equation}
Since the initial data $B(-\infty,y_2,\varphi)=
B_+(-\infty,y_2,\varphi)-B_-(-\infty,y_2,\varphi)$ are known we can
recover $B(x,e^{i\varphi})$ as the solution of the Cauchy problem.
The continuation of the proof is similar to [N] (see also \S 3).
Consider
\begin{equation}                                     \label{eq:4.25}
I(x)=\frac{1}{2\pi i}\int_{|t|=1}(B_+(x,t)-B_-(x,t))dt.
\end{equation}
We have $\int_{|t|=1}B_-(x,t)dt=0$ since $B_-(x,t)$ is analytic
when $|t|<1$ and continuous when $|t|\leq 1$.  It follows from
(\ref{eq:4.18}) that
\begin{equation}                                     \label{eq:4.26}
B_+(x,t)=c_+(x,\infty)\frac{1}{t}
\left(S(c_+^{-1}(x,\infty)V(x)c_+(x,\infty))\right)c_+^{-1}(x,\infty)+
O\left(\frac{1}{t^2}\right),
\end{equation}
where $t\rightarrow\infty,\ S$ is the same as in (\ref{eq:3.19}).
Therefore
\begin{equation}                                     \label{eq:4.27}
\frac{1}{2\pi i}\int_{|t|=1}B_+(x,t)dt
=c_+(x,\infty)
\left(S(c_+^{-1}(x,\infty)V(x)c_+(x,\infty))\right)c_+^{-1}(x,\infty).
\end{equation}
Multiplying (\ref{eq:4.27}) by $c_+^{-1}(x,\infty)$ from the left and 
by $c_+(x,\infty)$ from the right and apply the operator 
$\frac{\partial}{\partial\bar{z}}$ we get
\[
c_+^{-1}(x,\infty)V(x)c_+(x,\infty)
=\frac{\partial}{\partial\bar{z}}
\left(c_+^{-1}(x,\infty)I(x)c_+(x,\infty)\right)
\]
Finally
\[
V(x)=c_+(x,\infty)\left(
\frac{\partial}{\partial\bar{z}}
\left(c_+^{-1}(x,\infty)I(x)c_+(x,\infty)\right)\right)c_+^{-1}(x,\infty).
\]
\qed

\end{document}